\centerline{\bf Computational solutions of distributed order reaction-diffusion systems}
\vskip0cm\centerline{\bf associated with Riemann-Liouville derivatives}

\vskip.3cm\centerline{R.K. Saxena$^{d}$,~A.M. Mathai$^{b,c}$ and H.J. Haubold$^{a,b}$}
\vskip.3cm\noindent $^a$ Office for Outer Space Affairs, United Nations, P.O. Box 500, A-1400 Vienna, Austria,
\vskip.0cm\noindent hans.haubold@unvienna.org
\vskip0cm\noindent $^b$ Centre for Mathematical Sciences, Arunapuram P.O., Pala-686574, Kerala, India
\vskip0cm\noindent$^c$ Department of Mathematics and Statistics, McGill University, Montreal, Canada H3A 2K6,
\vskip.0cm\noindent mathai@math.mcgill.ca
\vskip0cm\noindent$^d$ Department of Mathematics and Statistics, Jai Narain Vyas University, Jodhpur-342004, India,
\vskip.0cm\noindent ram.saxena@yahoo.com

\vskip.5cm\noindent{\bf Abstract}

\vskip.3cm This article is in continuation of our earlier article [37] in which computational solution of an unified reaction-diffusion equation of distributed order associated with Caputo derivatives as the time-derivative and Riesz-Feller derivative as space derivative is derived. In this article, we present computational solutions of distributed order fractional reaction-diffusion equations associated with Riemann-Liouville derivatives of fractional orders as the time-derivatives and Riesz-Feller fractional derivatives as the space derivatives. The method followed in deriving the solution is that of joint Laplace and Fourier transforms. The solution is derived in a closed and computational form in terms of the familiar Mittag-Leffler function. It provides an elegant extension of the results given earlier by Chen et al. [1], Debnath [3], Saxena et al. [36], Haubold et al. [15] and Pagnini and Mainardi [30]. The results obtained are presented in the  form of two theorems. Some interesting results associated with fractional Riesz derivatives are also derived as special cases of our findings. It will be seen that in case of distributed order fractional reaction-diffusion, the solution comes in a compact and closed form in terms of a generalization of the Kamp\'e de F\'eriet hypergeometric series in two variables, defined by Srivastava and Daoust [46] (also see Appendix B). The convergence of the double series occurring in the solution is also given.



\vskip.5cm\noindent{\bf 1.\hskip.3cm Introduction}
\vskip.3cm Distributed order sub-diffusion is discussed by Naber [26]. Distributed order fractional diffusion systems are studied, among others, by Saxton [40,41,42], Langlands [21], Sokolov et al. [43,44], Sokolov and Klafter [44], Saxena and Pagnini [38] and Nikolova and Boyadijiev [28], and recent monographs on the subject [4,20,22,25]. General models for reaction-diffusion systems are discussed by Wilhelmsson and Lazzaro [46], Henry and Wearne [17,18]. Henry et al. [19], Mainardi et al. [23,24], Haubold et al. [14,15,16], Saxena et al. [34,35,36,37,38] and others. Stability in reaction-diffusion systems and nonlinear oscillations have been discussed by Gafiychuk et al. [10,11]. Pattern formation in reaction-diffusion as well as non-Gaussian, non-Markovian, and non-Fickian phenomena related to astronomical, physical, chemical, and biological sciences can be found in the work of Cross and Hohenberg [2], Nicolis and Prigogine [27], and Haubold et al. [16]. Recently, Engler [6] discussed the speed of spread of fractional reaction-diffusion.

\vskip.2cm In a recent article, Chen et al. [1] have derived the fundamental and numerical solution of a reaction-diffusion equation associated with the Riesz fractional derivative as the space derivative. Reaction-diffusion models associated with Riemann-Liouville fractional derivative as the time derivative and Riesz-Feller derivative as the space derivative are recently discussed by Haubold et al. [14]. Such equations in case of Caputo fractional derivative are solved by Saxena et al. [37].
\vskip.2cm The main object of this article is to investigate the computational solutions of fraction reaction-diffusion equations (2.1) and (4.1). The results are obtained in a closed and computational forms. Due to general character of the derived results, many known results given earlier by Chen et al. [1], Haubold et al. [15] and Pagnini and Mainardi [30], Saxena et al. [37], and others, readily follow as special cases of our derived results.

\vskip.3cm\noindent{\bf 2.\hskip.3cm Solution of fractional reaction-diffusion equations}
\vskip.3cm In this section, we will investigate a computational solution of  one-dimensional fractional reaction-diffusion equation (2.1) containing Riemann-Liouville derivative as the time-derivative and a finite number of Riesz-Feller derivative as the space derivatives. The result obtained are in a compact and computational form in terms of the generalized Mittag-Leffler function, defined by (2.5) in the form of the following theorem:

\vskip.3cm\noindent{\bf Theorem 1.}\hskip.3cm{\it Consider the following one-dimensional non-homogeneous unified fractional reaction-diffusion model associated with time-derivative ${_0D_t^{\alpha}}$ defined by (A1), $n\in N$, and Riesz-Feller space-derivatives ${_xD_{\theta_1}^{\gamma_1}},...,{_xD_{\theta_n}^{\gamma_n}}$, defined by (A3):

$$\eqalignno{{_0D_{t}^{\alpha}}N(x,t)&=\sum_{j=1}^n\mu_j~{_xD_{\theta_j}^{\gamma_j}}N(x,t)+\phi(x,t)&(2.1)\cr
\noalign{\hbox{where $t>0,x\in R; \alpha,\theta_1,...,\theta_n,\gamma_1,...,\gamma_n$ are real parameters with the condition}}
\mu_j>0,0<\gamma_j&\le 2, j=1,...,n,~|\theta_j|\le \min_{1\le j\le n}(\gamma_j,2-\gamma_j), 1<\alpha\le 2;&(2.2)\cr
\noalign{\hbox{with initial conditions}}
{_0D_t^{\alpha-1}}N(x,0)&=f(x),~{_0D_t^{\alpha-2}}N(x,0)=g(x),~x\in R,~\lim_{x\rightarrow\pm\infty}N(x,t)=0,t>0.&(2.3)\cr}
$$Here ${_0D_t^{\nu}}N(x,t)$ means the Riemann-Liouville fractional partial derivative of $N(x,t)$ with respect to $t$ of order $\nu$ evaluated at $t=0,\nu=\alpha-1,\alpha-2$; ${_xD_{\theta_1}^{\gamma_1}},...,{_xD_{\theta_n}^{\gamma_n}}$ are the Riesz-Feller space fractional derivatives with asymmetries $\theta_1,...,\theta_n$ respectively. Further, ${_0D_t^{\alpha}}$ is the Riemann-Liouville time-fractional derivative of order $\alpha$; $\mu_1,...,\mu_n$ are arbitrary constants; $f(x), g(x)$ and $\phi(x,t)$ are given functions. Then for the solution of (2.1), subject to the above conditions, there holds the formula
$$\eqalignno{N(x,t)&={{f^{\alpha-1}}\over{2\pi}}\int_{-\infty}^{\infty}t^{*}(k)
E_{\alpha,\alpha}(-bt^{\alpha})\exp(-ikx){\rm d}k\cr
&+{{t^{\alpha-2}}\over{2\pi}}\int_{-\infty}^{\infty}g^{*}(k)E_{\alpha,\alpha-1}(-bt^{\alpha})\exp(-ikx){\rm d}k\cr
&+{{1}\over{2\pi}}\int_0^t\xi^{\alpha-1}\int_{-\infty}^{\infty}\phi^{*}(k,t-\xi)E_{\alpha,\alpha}(-b\xi^{\alpha})
\exp(-ikx){\rm d}k~{\rm d}\xi,&(2.4)\cr}
$$where $E_{\alpha,\beta}(z)$ is the generalized Mittag-Leffler function, defined by Wiman [48] (also see Erd\'elyi et al. [7], Dzherbashyan [5]) in the following form:
$$\eqalignno{E_{\alpha,\beta}(z)&=\sum_{n=0}^{\infty}{{z^n}\over{\Gamma(\alpha n+\beta)}},~\alpha,\beta\in C, \Re(\alpha)>0,\Re(\beta)>0&(2.5)\cr
\noalign{\hbox{and}}
b&=\sum_{j=1}^n\mu_j~\psi_{\gamma_j}^{\theta_j}(k).&(2.6)\cr}
$$}

\vskip.3cm\noindent{\bf Proof:}\hskip.3cm In order to derive the solution of (2.1), we introduce the joint Laplace-Fourier transform int he form
$${\tilde{N}}^{*}(k,s)=\int_0^{\infty}\int_{-\infty}^{\infty}{\rm e}^{-st+ikx}N(x,t){\rm d}~{\rm d}t
$$where $\Re(s)>0,k>0$. If we apply the Laplace transform with respect to the time variable $t$, Fourier transform with respect to space variable $x$, use the initial conditions (2.2), (2.3) and the formula (A2), then the given equation transforms into the form
$$s^{\alpha}{\tilde{N}}^{*}(k,s)-f^{*}(k)-sg^{*}(k)=-\sum_{j=1}^n\mu_j~\psi_{\gamma_j}^{\theta_j}(k){\tilde{N}}^{*}(k,s)
+{\tilde{\phi}}^{*}(k,s)
$$where according to the conventions followed, the symbols $\sim$ will stand for the Laplace transform with respect to time variable $t$ and $*$ represents the Fourier transform with respect to space variable $x$. Solving for ${\tilde{N}}^{*}(k,s)$ it yields
$$\eqalignno{{\tilde{N}}^{*}(k,s)&={{f^{*}(k)}\over{s^{\alpha}+b}}+{{sg^{*}(k)}\over{s^{\alpha}+b}}
+{{{\tilde{\phi}}^{*}(k,s)}\over{s^{\alpha}+b}}&(2.7)\cr
\noalign{\hbox{where}}
b&=\sum_{j=1}^n\mu_j\psi_{\gamma_j}^{\theta_j}(k).\cr}
$$On taking the inverse Laplace transform of (2.7) by means of the formula

$$L^{-1}\left[{{s^{\alpha-1}}\over{b+s^{\beta}}};t\right]=t^{\beta-\alpha}E_{\beta,\beta-\alpha+1}(-bt^{\beta})\eqno(2.8)
$$where $\Re(s)>0,\Re(\beta)>0,\Re(\beta-\alpha)>-1;|{{b}\over{s^{\beta}}}|<1$ it is found that
$$\eqalignno{N^{*}(k,t)&=t^{\alpha-1}f^{*}(k)E_{\alpha,\alpha}(-bt^{\alpha})+t^{\alpha-2}g^{*}(k)
E_{\alpha,\alpha-1}(-bt^{\alpha})\cr
&+\int_0^t\phi^{*}(k,t-\xi)\xi^{\alpha-1}E_{\alpha,\alpha}(-b\xi^{\alpha}){\rm d}\xi.&(2.9)\cr}
$$The required solution (2.4) is now obtained by taking the inverse Fourier transform of (2.9). This completes the proof of Theorem 1.

\vskip.3cm\noindent{\bf 3.\hskip.3cm Special cases}
\vskip.3cm\noindent (i): ~If we take $\theta_1=...=\theta_n=0$ then by virtue of the relation (A11), we obtain the following corollary.

\vskip.3cm\noindent{\bf Corollary 1.1.}\hskip.3cm{\it Consider the following one-dimensional non-homogeneous unified fractional reaction-diffusion model associated with time derivative ${_0D_t^{\alpha}}$ obtained by (A1), $n\in N$, and Riesz space derivatives ${_xD_0^{\gamma_1}},...,{_xD_0^{\gamma_n}}$ defined by (A3):
$${_0D_t^{\alpha}}N(x,t)=\sum_{j=1}^n\mu_j~{_xD_0^{\gamma_j}}N(x,t)+\phi(x,t),\eqno(3.1)
$$where $t>0,x\in R; \alpha,\gamma_1,...,\gamma_n$ are real parameters with the constraints
$$\mu_j>0,~0<\gamma_j\le 2, j=1,...,n, 1<\alpha\le 2;\eqno(3.2)
$$with initial conditions
$${_0D_t^{\alpha-1}}N(x,0)=f(x),~{_0D_t^{\alpha-2}}N(x,0)=g(x),~x\in R,\lim_{x\rightarrow\pm\infty}N(x,t)=0,t>0\eqno(3.3)
$$where the various terms are as defined in (2.3). Then there holds the formula
$$\eqalignno{N(x,t)&={{t^{\alpha-1}}\over{2\pi}}\int_{-\infty}^{\infty}f^{*}(k)E_{\alpha,\alpha}(-ct^{\alpha})
\exp(-ikx){\rm d}k\cr
&+{{t^{\alpha-2}}\over{2\pi}}\int_{-\infty}^{\infty}g^{*}(k)E_{\alpha,\alpha-1}(-ct^{\alpha})\exp(-ikx){\rm d}k\cr
&+{{1}\over{2\pi}}\int_0^t\xi^{\alpha-1}\int_{-\infty}^{\infty}\phi^{*}(k,t-\xi)E_{\alpha,\alpha}(-ct^{\alpha})
\exp(-ikx){\rm d}k~{\rm d}\xi&(3.4)\cr}
$$where $E_{\alpha,\beta}(z)$ is the generalized Mittag-Leffler function, defined by (2.5) and $c=\sum_{j=1}^n\mu_j|k|^{\gamma_j}$.}

\vskip.3cm\noindent Note that when $g(x)=0$ then Theorem 1 yields (2.9) where the middle term involving $g^{*}(k)$ will be absent.

\vskip.3cm\noindent{\bf Corollary 1.2.}\hskip.3cm{\it Consider the same equation in (3.1) under the conditions
$$\eqalignno{\mu_j>0,~0<\gamma_j\le 2,j=1,...,n,&|\theta_j|\le \min_{1\le j\le n}(\gamma_j,2-\gamma_j),~0<\alpha\le 1;&(3.6)\cr
\noalign{\hbox{with the initial conditions}}
{_0D_t^{\alpha-1}}N(x,0)&=\delta(x),~{_0D_t^{\alpha-2}}N(x,0)=0,~x\in R,\lim_{x\rightarrow\pm\infty}N(x,t)=0,t>0&(3.7)\cr}
$$where $\delta(x)$ is Dirac delta function. Then for the fundamental solution there holds the formula

$$\eqalignno{N(x,t)&={{t^{\alpha-1}}\over{2\pi}}\int_{-\infty}^{\infty}E_{\alpha,\alpha}(-bt^{\alpha})\exp(-ikx){\rm d}k\cr
&+{{1}\over{2\pi}}\int_0^t\xi^{\alpha-1}\int_{-\infty}^{\infty}\phi^{*}(k,t-\xi)E_{\alpha,\alpha}(-b\xi^{\alpha})
\exp(-ikx){\rm d}k~{\rm d}\xi&(3.8)\cr}
$$where $E_{\alpha,\beta}(z)$ is the generalized Mittag-Leffler function defined by (2.5) and $b$ is given in (2.6).}

\vskip.3cm If we further set $\phi=0,n=1,\mu_1=\mu,\theta_1=\theta,\gamma_1=\gamma$ and apply the formula (A14) then we have the following result:

\vskip.3cm\noindent{\bf Corollary 1.3.}\hskip.3cm{\it Consider the following reaction-diffusion model
$$\eqalignno{{_0D_t^{\alpha}}N(x,t)&=\mu~{_xD_{\theta}^{\gamma}}N(x,t)&(3.9)\cr
\noalign{\hbox{with the initial conditions}}
{_0D_t^{\alpha-1}}N(x,0)&=\delta(x),~x\in R,~0<\alpha\le 1,\lim_{x\rightarrow\pm\infty}N(x,t)=0,t>0&(3.10)\cr}
$$where $\mu$ is a diffusion constant, $\mu,t>0; \alpha,\beta,\theta$ are real parameters with the constraint
$$0<\gamma\le 2, ~|\theta|\le \min (\gamma,2-\gamma)\eqno(3.11)
$$and $\delta(x)$ is a Dirac delta function. Then the fundamental solution of (3.9) with initial conditions, there holds the formula
$$N(x,t)={{t^{\alpha-1}}\over{\gamma|x|}}H_{3,3}^{2,1}\left[{{|x|}\over{(\mu t^{\alpha})^{1/\gamma}}}\bigg\vert_{(1,1/\gamma), (1,1), (1,\rho)}^{(1,1/\gamma),(\alpha,\alpha/\gamma), (1,\rho)}\right],~\gamma>0\eqno(3.12)
$$where $H_{3,3}^{2,1}(\cdot)$ is the familiar H-function defined by (A12), also see [25]; $\rho={{\alpha-\theta}\over{2\alpha}}$.}

\vskip.3cm For $\theta=0$, (3.12) reduces to the result given by Saxena et al. [38] in a slightly different form. On the other hand, if we further set $\alpha=1$ we obtain the following result given by Chen et al. [1] in a different form:

\vskip.3cm\noindent{\bf Corollary 1.4.}\hskip.3cm{\it Consider the following reaction-diffusion model
$$\eqalignno{{{{\rm d}}\over{{\rm d}t}}N(x,t)&=\mu~{_xD_0^{\gamma}}N(x,t)&(3.13)\cr
\noalign{\hbox{with the initial consitions}}
N(x,0)&=\delta(x), x\in R,\lim_{x\rightarrow\pm\infty}N(x,t)=0,t>0,0<\beta\le 2,&(3.14)\cr}
$$where $\mu$ is a diffusion constant, $\mu,t>0,\beta$ are real parameters and $\delta(x)$ is Dirac delta function. Then for the fundamental solution of (3.13) with initial conditions (3.14) there holds the formula
$$N(x,t)={{1}\over{\gamma|x|}}H_{2,2}^{1,1}\left[{{|x|}\over{(\mu t)^{1/\gamma}}}\bigg\vert_{(1,1),(1,1/2)}^{(1,1/\gamma),(1,1/2)}\right],\gamma>0.\eqno(3.15)
$$}

\vskip.3cm\noindent{\bf Remark 3.1.}\hskip.3cm The result obtained by Chen et al. [1] is in terms of the Fourier transform of the Mittag-Leffler function, whereas our result (3.15) is in terms of the H-function in a closed and computable form. It is interesting to observe that for $n=1$, Theorem 1 reduces to one given by Haubold et al. [15].

\vskip.3cm\noindent{\bf 4.\hskip.3cm Further results on distributed order reaction-diffusion systems}

\vskip.3cm In this section we will investigate a computational solution of a distributed order reaction-diffusion equation containing two Riemann-Liouville derivatives. The solution is obtained in terms of generalized Mittag-Leffler function due to Prabhakar [32]. The main result can be expressed in terms of Srivastava-Daoust hypergeometric function of two variables [45].

\vskip.3cm\noindent{\bf Theorem 2.}\hskip.3cm{\it Consider the following unified one-dimensional non-homogeneous reaction-diffusion equation of fractional order:
$${_0D_t^{\alpha}}N(x,t)+a~{_0D_t^{\beta}}N(x,t)=\sum_{j=1}^n\mu_j~{_xD_{\theta_j}^{\gamma_j}}N(x,t)+\phi(x,t)\eqno(4.1)
$$where $t>0,a,x\in R; \alpha,\beta,\theta_1,...,\theta_n$ are real parameters with the constraints
$$\mu_j>0,~0<\gamma_j\le 2, j=1,...,n,~|\theta_j|\le \min_{1\le j\le n}(\gamma_j,2-\gamma_j),~1<\alpha\le 2,~1<\beta\le 2\eqno(4.2)
$$with the initial conditions
$$\eqalignno{{_0D_t^{\alpha-1}}N(x,0)&=f_1(x),~{_0D_t^{\alpha-2}}N(x,0)=g_1(x),~{_0D_t^{\beta-1}}N(x,0)=f_2(x),\cr
{_0D_t^{\beta-2}}N(x,0)&=g_2(x),~x\in R,\lim_{x\rightarrow\pm\infty}N(x,t)=0,t>0&(4.3)\cr}
$$where the various quantities are as defined in Theorem 1, and $f_1(x),f_2(x),g_1(x),g_2(x)$ and $\phi(x,t)$ are given functions. Then for the solution of (4.1), subject to the above conditions (4.2),(4.3), there holds the formula
$$\eqalignno{N(x,t)&={{t^{\alpha-1}}\over{2\pi}}\int_{-\infty}^{\infty}(f_1^{*}(k)+af_2^{*}(k))
\sum_{r=0}^{\infty}(-a)^rt^{(\alpha-\beta)r}E_{\alpha,\alpha+(\alpha-\beta)r}^{r+1}(-bt^{\alpha})\exp(-ikx){\rm d}k\cr
&+{{t^{\alpha-2}}\over{2\pi}}\int_{-\infty}^{\infty}(g_1^{*}(k)+ag_2^{*}(k))\sum_{r=0}^{\infty}(-a)^rt^{(\alpha-\beta)r}
E_{\alpha,\alpha+(\alpha-\beta)r-1}^{r+1}(-bt^{\alpha})\exp(-ikx){\rm d}k\cr
&+{{1}\over{2\pi}}\int_0^t\xi^{\alpha-1}\int_{-\infty}^{\infty}\sum_{r=0}^{\infty}(-a)^ru^{(\alpha-\beta)r}
\int_{-\infty}^{\infty}\phi^{*}(p,t-\xi)E_{\alpha,\alpha+(\alpha-\beta)r}^{r+1}(-b\xi^{\alpha})\exp(-ikx){\rm d}k~{\rm d}\xi&(4.4)\cr}
$$where $E_{\alpha,\beta}^{\gamma}(z)$ is the generalized Mittag-Leffler function ([32])
$$E_{\alpha,\beta}^{\gamma}(z)=\sum_{n=0}^{\infty}{{(\gamma)_nz^n}\over{\Gamma(\alpha n+\beta)n!}},~\Re(\alpha)>0,\Re(\beta)>0\eqno(4.5)
$$where $b=\sum_{j=1}^n\mu_j\psi_{\gamma_j}^{\theta_j}(k)$.}

\vskip.3cm\noindent{\bf Proof:}\hskip.3cm If we apply the Laplace transform with respect to the time variable $t$, Fourier transform with respect to space variable $x$ and use the initial conditions (4.2),(4.3) and the formula (4.5), then the given equation transforms into the form
$$\eqalignno{s^{\alpha}{\tilde{N}}^{*}(k,s)&-f_1^{*}(k)-sg_1^{*}(k)+as^{\beta}{\tilde{N}}^{*}(k,s)-af_2^{*}(k)-asg_2^{*}(k)\cr
&=-\sum_{j=1}^n\mu_j\psi_{\gamma_j}^{\theta_j}(k){\tilde{N}}^{*}(k,s)+{\tilde{\phi}}^{*}(k,s).&(4.6)\cr
\noalign{\hbox{Solving for ${\tilde{N}}^{*}(k,s)$ it yields}}
{\tilde{N}}^{*}(k,s)&={{f_1^{*}(k)+af_2^{*}(k)}\over{s^{\alpha}+as^{\beta}+b}}+{{s(g_1^{*}(k)+ag_2^{*}(k))}
\over{s^{\alpha}
+as^{\beta}+b}}+{{{\tilde{\phi}}^{*}(k)}\over{s^{\alpha}+as^{\beta}+b}}&(4.7)\cr}
$$where $b$ is defined in (2.6). To invert the equation (4.7), we first invert the Laplace transform and then the Fourier transform. Thus to invert the Laplace transform we use the formula given in [40]:
$$L^{-1}\left[{{s^{\rho-1}}\over{s^{\alpha}+as^{\beta}+b}}\right]=t^{\alpha-\rho}\sum_{r=0}^{\infty}(-a)^r
t^{(\alpha-\beta)r}E_{\alpha,\alpha+(\alpha-\beta)r-\rho+1}^{r+1}(-bt^{\alpha})\eqno(4.8)
$$where $\Re(s)>0,\Re(\alpha)>0,\Re(\beta)>0,\Re(\alpha-\rho)>-1,\Re(\alpha-\beta)>0$, $|{{as^{\beta}}\over{b+s^{\alpha}}}|<1$; and the convolution theorem of the Laplace transform to obtain
$$\eqalignno{N^{*}(k,t)&=t^{\alpha-1}[f_1^{*}(k)+af_2^{*}(k)]\sum_{r=0}^{\infty}(-a)^rt^{(\alpha-\beta)r}
E_{\alpha,\alpha+(\alpha-\beta)r}^{r+1}(-bt^{\alpha})\cr
&+t^{\alpha-2}[g_1^{*}(k)+ag_2^{*}(k)]\sum_{r=0}^{\infty}(-a)^rt^{(\alpha-\beta)r}
E_{\alpha,\alpha+(\alpha-\beta)r-1}^{r+1}(-bt^{\alpha})\cr
&+\sum_{r=0}^{\infty}(-a)^r\int_0^t\phi^{*}(k,t-\xi)\xi^{\alpha+(\alpha-\beta)r-1}
E_{\alpha,\alpha+(\alpha-\beta)r}^{r+1}(-b\xi^{\alpha}){\rm d}\xi.&(4.9)\cr}
$$Now, the application of the inverse Fourier transform gives the following solution:

$$\eqalignno{N(x,t)&={{t^{\alpha-1}}\over{2\pi}}\int_{-\infty}^{\infty}[f_1^{*}(k)+ag_2^{*}(k)]
\sum_{r=0}^{\infty}(-a)^rt^{(\alpha-\beta)r}E_{\alpha,\alpha+(\alpha-\beta)r}^{r+1}(-bt^{\alpha})\exp(-ikx){\rm d}k\cr
&+{{t^{\alpha-2}}\over{2\pi}}\int_{-\infty}^{\infty}[g_1^{*}(k)+ag_2^{*}(k)]\sum_{r=0}^{\infty}(-a)^rt^{(\alpha-\beta)r}
E_{\alpha,\alpha+(\alpha-\beta)r-1}^{r+1}(-bt^{\alpha})\exp(-ikx){\rm d}k\cr
&+\sum_{r=0}^{\infty}{{(-a)^r}\over{2\pi}}\int_0^t\xi^{\alpha+(\alpha-\beta)r-1}\int_{-\infty}^{\infty}\phi^{*}(k,t-\xi)
\exp(-ikx)\cr
&\times E_{\alpha,\alpha+(\alpha-\beta)r}^{r+1}(-b\xi^{\alpha}){\rm d}k.&(4.10)\cr}
$$

\vskip.3cm\noindent{\bf Alternative form of the solution (4.4)}

\vskip.3cm By virtue of the series representation of the generalized Mittag-Leffler function $E_{\beta,\gamma}^{\alpha}(z)$ defined in (4.5), the expression
$$\eqalignno{t^{\alpha-\rho}&\sum_{r=0}^{\infty}(-a)^rt^{(\alpha-\beta)r}
E_{\alpha,\alpha+(\alpha-\beta)r-\rho+1}^{r+1}(-bt^{\alpha})&(4.11)\cr
\noalign{\hbox{can be written as}}
t^{\alpha-\rho}&\sum_{r=0}^{\infty}\sum_{u=0}^{\infty}{{(1)_{r+u}}\over{r!u!}}{{(-at^{\alpha-\beta})^r(-bt^{\alpha})^u}
\over{\Gamma(\alpha-\rho+1+(\alpha-\beta)r+au)}}\cr
&=t^{\alpha-\rho}S_{1:0;0}^{1:0;0}\left[-at^{\alpha-\beta},-bt^{\alpha}\bigg\vert_{[\alpha-\rho+1:\alpha-\beta;\alpha]:-;-}^{[1:1;1]:-;-}\right],&(4.12)\cr}
$$where $S(\cdot)$ is the Srivastava-Daoust generalization of the Kamp\'e de F\'eriet hypergeometric series in two variables [46], for its definition, see (B2) in Appendix B. Hence, Theorem 2 can now be stated in terms of the Srivastava-Daoust hypergeometric function of two variables in the following form: Under the conditions of Theorem 2, the unified one-dimensional fractional reaction-diffusion equation

$$\eqalignno{{_0D_t^{\alpha}}N(x,t)&+a~{_0D_t^{\beta}}N(x,t)=\lambda~{_xD_{\theta}^{\gamma}}N(x,t)+\phi(x,t),&(4.13)\cr
\noalign{\hbox{has the solution given by}}
N(x,t)&={{t^{\alpha-1}}\over{2\pi}}\int_{-\infty}^{\infty}[f_1^{*}(k)+af_2^{*}(k)]\exp(-ikx)\cr
&\times S_{1:0;0}^{1:0;0}\left[-at^{\alpha-\beta},-bt^{\alpha}\bigg\vert_{[\alpha:\alpha-\beta;\alpha]:-:-}^{[1:1;1]:-:-}\right]
{\rm d}k\cr
&+{{t^{\alpha-2}}\over{2\pi}}\int_{-\infty}^{\infty} [g_1^{*}(k)+ag_2^{*}(k)]S_{1:0;0}^{1:0;0}\left[-at^{\alpha-\beta},-bt^{\alpha}\bigg\vert_{[\alpha-1:\alpha-\beta;\alpha]
:-;-}^{[1:1;1]:-;-}\right]{\rm d}k\cr
&+{{1}\over{2\pi}}\int_0^t\xi^{\alpha-1}\int_{-\infty}^{\infty}\phi^{*}(k,t-\xi)\exp(-ikx)\cr
&\times S_{1:0;0}^{1:0;0}
\left[-au^{\alpha-\beta},-b\xi^{\alpha}\bigg\vert_{[\alpha:\alpha-\beta;\alpha]:-;-}^{[1:1;1]-;-}\right]{\rm d}k~{\rm d}\xi&(4.14)\cr}
$$where $\Re(\alpha)>0,\Re(\beta)>0,\Re(\alpha-\beta)>0$ and $b$ is defined in (2.6).

\vskip.3cm\noindent{\bf Note 4.1.}\hskip.3cm By virtue of the lemma given in the Appendix B, the double infinite power series occurring in Theorem 2 converges for $\Re(\alpha)>0,\Re(\alpha-\beta)>0$.

\vskip.3cm\noindent{\bf 5.\hskip.3cm Special cases of Theorem 2}

\vskip.3cm If we set $g_1(x)=g_2(x)=0$ in Theorem 2 then it reduces to the following:

\vskip.3cm\noindent{\bf Corollary 2.1.}.\hskip.3cm{\it Consider the following unified one-dimensional non-homogeneous reaction-diffusion equation of fractional order:
$${_0D_t^{\alpha}}N(x,t)+a~{_0D_t^{\beta}}N(x,t)=\sum_{j=1}^n\mu_j~{_xD_{\theta_j}^{\gamma_j}}N(x,t)+\phi(x,t)\eqno(5.1)
$$where $t>0,a,x\in R; \alpha,\beta,\theta_1,...,\theta_n, \gamma_1,...,\gamma_n$ are real parameters with the constraint
$$\mu_j>0,0<\gamma_j\le 2, j=1,...,n,|\theta_j|\le \min_{1\le j\le n}(\gamma_j,2-\gamma_j), 1<\alpha\le 2, 1<\beta\le 2\eqno(5.2)
$$with the initial conditions
$$\eqalignno{{_0D_t^{\alpha-1}}N(x,0)&=f_1(x),~{_0D_t^{\alpha-2}}N(x,t)=0,~{_0D_t^{\beta-1}}N(x,0)=f_2(x),\cr
{_0D_t^{\beta-2}}N(x,0)&=0,~x\in R, \lim_{x\rightarrow\pm\infty}N(x,t)=0,t>0&(5.3)\cr}
$$where the various quantities are as defined in Theorem 2. Then we have a special case of (4.10) and (4.12) with the corresponding changes.}

\vskip.3cm Other special cases of interest are the situations (i): $n=1$, (ii): $\theta_1=...=\theta_n=0$, (iii): $f_1(x)=f_2(x)=\delta(x)$  in Theorem 2. We get the corresponding results from (4.10) and (4.14) by substitution.

\vskip.5cm\noindent{\bf References}

\vskip.3cm\noindent [1] Chen, J., Liu, F., Turner, I. and Anh, V. : The fundamental and numerical solutions of the Riesz space-fractional reaction-dispersion equation, {\it The Australian and New Zealand Industrial and Applied Mathematics Journal (ANZIAM)}, {\bf 50}(2008), 45-57.

\vskip.3cm\noindent [2] Cross, M.C. and Hohenberg, P.C. : Pattern formation outside of equilibrium, {\it Reviews of Modern Physics}, {\bf 65}(1993), 851-912.

\vskip.3cm\noindent [3] Debnath, L. : Fractional integral and fractional differential equations in fluid mechanics, {\it Frac. Calc. Appl. Anal.}, {\bf 6}(2003), 119-155.

\vskip.3cm\noindent [4] Diethelm, K.: {\it The Analysis of Fractional Differential Equations}, Springer, Berlin, Heidelberg, 2010.

\vskip.3cm\noindent [5] Dzherbashyan, M.M.: {\it Harmonic Analysis and Boundary Value Problems in the Complex Domain}, Birkha\"user-Verlag, Basel and London, 1993.

\vskip.3cm\noindent [6] Engler, H.: On the speed of spread for fractional reaction-diffusion, {\it International Journal of Differential Equations}, {\bf 2010 Article ID 315421}, 16 pages.

\vskip.3cm\noindent [7] Erd\'elyi, A., Magnus, W., Oberhettinger, F. and Tricomi, F.G. : {\it Higher Transcendental Functions}, {\bf Vol.3}, McGraw-Hill, New York, 1955.

\vskip.3cm\noindent [8] Feller, W. : On a generalization of Marcel Riesz' potentials and the semi-groups generated by them, {\it Meddelanden Lunds Universitets Matematiska Seminarium (Comm. S\'em. Math\'em. Universit\'e de Lund)}, {\bf Tome Suppl. D\'edi\'e \'a M. Riesz, Lund}, 73-81(1952).

\vskip.3cm\noindent [9] Feller, W.: {\it An Introduction to Probability Theory and Its Applications}, {\bf Vol.2}, Second Edition, New York, Wiley, 1971.

\vskip.3cm\noindent [10] Gafiychuk, V., Datsko, B. and Meleshko, V. : Mathematical modeling in pattern formation in sub and super-diffusive reaction-diffusion systems, 2006, arXiv:nlinAO/0811005v3.

\vskip.3cm\noindent [11] Gafiychuk, V., Datsko, B. and Meleshko, V. : Nonlinear oscillations and stability domains in fractional reaction-diffusion systems, 2007, arXiv:nlinPS/0702013v1.

\vskip.3cm\noindent [12] Gorenflo, R. and Mainardi, F. : Approximation of Levy-Feller diffusion by random walk, {\it Journal for Analysis and Its Applications}, {\bf 18(2)}(1999), 1-16.

\vskip.3cm\noindent [13] Guo, X, and Xu, M. : Some physical applications of Schr\"odinger equation, {\it J. Math. Phys.}, {\bf 47082104}(2008):doi10, 1063/1.2235026 (9 pages).

\vskip.3cm\noindent [14] Haubold , H.J., Mathai, A.M. and Saxena, R.K. : Solutions of reaction-diffusion equations in terms of the H-function, {\it Bull. Astro. Soc. India}, {\bf 35(4)}(2007), 681-689.

\vskip.3cm\noindent [15] Haubold, H.J., Mathai, A.M. and Saxena, R.K.: Further solutions of reaction-diffusion equations in terms of the H-function, {\it J. Comput. Appl. Math.}, {\bf 235}(2011), 1311-1316.

\vskip.3cm\noindent [16] Haubold, H.J., Mathai, A.M. and Saxena, R.K. : Analysis of solar neutrino data from SuperKamiokande I and II: Back to the solar neutrino problem, 2012, arXiv:astro-ph.SR/1209.1520.

\vskip.3cm\noindent [17] Henry, B.I. and Wearne, S.: Fractional reaction-diffusion, {\it Physica A}, {\bf 276}(2000), 448-455.

\vskip.3cm\noindent [18] Henry, B.I. and Wearne, S. : Existence of Turing instabilities in a two species reaction-diffusion system, {\it SIAM Journal of Applied Mathematics}, {\bf 62}(2002), 870-887.

\vskip.3cm\noindent [19] Henry, B.I., Langlands, T.A.M. and Wearne, S.L. : Turing pattern formation in fractional activator-inhibitor systems, {\it Phys. Rev.},{\bf E72}(2005), 026101. 

\vskip.3cm\noindent [20] Kilbas, A.A., Srivastava, H.M. and Trujillo, J.J.:{\it Theory and Applications of Fractional Differential Equations}, Elsevier, Amsterdam, 2006.

\vskip.3cm\noindent [21] Langlands, T.A.M.: Solution of a modified fractional diffusion equation, {\it Physica A}, {\bf 367}(2006), 136-144.

\vskip.3cm\noindent [22] Mainardi, F. : {\it Fractional Calculus and Waves in Linear Viscoelasticity}, World Scientific, Singapore and Imperial College Press, 2010.

\vskip.3cm\noindent [23] Mainardi, F., Luchko, Y. and Pagnini, G. : The fundamental solution of the space-time fractional diffusion equation, {\it Fract. Calc. Appl. Anal.}, {\bf 4(2)}(2001), 153-202.

\vskip.3cm\noindent [24] Mainardi, F., Pagnini, G. and Saxena, R.K. : Fox H-functions in fractional diffusion, {\it J. Comput. Appl. Math.}, {\bf 178}(2005), 321-331.

\vskip.3cm\noindent [25] Mathai, A.M., Saxena, R.K. and Haubold, H.J.:{\it The H-function: Theory and Applications}, Springer, New York, 2010.

\vskip.3cm\noindent [26] Naber, M.: Distributed order fractional sub-diffusion, {\it Fractals}, {\bf 12(1)}(2004), 23-32.

\vskip.3cm\noindent [27] Nicolis, G. and Prigogine, I.: {\it Self-organization in Nonequilibrium Systems: From Dissipative Structures to Order Through Fluctuations}, Wiley, New York, 1977.

\vskip.3cm\noindent [28] Nikolova, Y. and Boyadijiev, L. : Integral transform methods to solve a time-space fractional diffusion equation, {\it Fract. Calc. Appl. Anal.}, {\bf 13}(2010), 57-67.

\vskip.3cm\noindent [29] Oldham, K.B. and Spanier, J.: {\it The Fractional Calculus Theory and Applications of Differentiation and Integration to Arbitrary Order}, Academic Press, New York, 1974.

\vskip.3cm\noindent [30] Pagnini, R. and Mainardi, F. : Evolution equations for a probabilistic generalization of Voigt profile function, {\it J. Comput. Appl. Math.}, {\bf 233}(2010), 1590-1595.

\vskip.3cm\noindent [31] Podlubny, I. : {\it Fractional Differential Equations}, Academic Press, New York, 1999.

\vskip.3cm\noindent [32] Prabhakar, T.R.: A singular integral equation with a generalized Mittag-Leffler function in the kernel, {\it Yokohama Math.J.}, {\bf 19}(1971), 7-15.

\vskip.3cm\noindent [33] Samko, S.G., Kilbas, A.A. and Marichev, O.I. : {\it Fractional Integrals and Derivatives: Theory and Applications}, Gordon and Breach Science Publishing, Switzerland, 1993.

\vskip.3cm\noindent [34] Saxena, R.K., Mathai, A.M. and Haubold, H.J. : Fractional reaction-diffusion equations, {\it Astrophysics and Space Science}, {\bf 305}(2006a), 289-296.

\vskip.3cm\noindent [35] Saxena, R.K., Mathai, A.M. and Haubold, H.J. : Reaction-diffusion systems and nonlinear waves, {\it Astrophysics and Space Science}, {\bf 305}(2006b), 297-303.

\vskip.3cm\noindent [36] Saxena, R.K., Mathai, A.M. and Haubold, H.J. : Solution of generalized reaction-diffusion equations, {\it Astrophysics and space Science}, {\bf 305}(2006c), 305-313.

\vskip.3cm\noindent [37] Saxena, R.K., Mathai, A.M. and Haubold, H.J. : Distributed order reaction-diffusion systems associated with Caputo derivatives, {\it arXiv:math-ph/1109.4841}.

\vskip.3cm\noindent [38] Saxena, R.K. and Pagnini, G. : Exact solutions of triple order time-fractional differential equations for anomalous relaxation and diffusion: The accelerating case, {\it Physica A}, {\bf 390}(2011), 602-613.

\vskip.3cm\noindent [39] Saxena, R.K., Saxena, R. and Kalla, S.L. : Computational solution of a fractional generalization of Schr\"odinger equation occurring in Quantum Mechanics, {\it Appl. Math. Comput.}, {\bf 216}(2010), 1412-1417.

\vskip.3cm\noindent [40] Saxton, M.: Anomalous diffusion due to obstacles: a Monte Carlo Study, {\it Biophys. J.}, {\bf 66}(1994), 394-401.

\vskip.3cm\noindent [41] Saxton, M.: Anomalous diffusion due to binding: a Monte Carlo Study, {\it Biophys. J.}, {\bf 70}(1996), 1250-1262.

\vskip.3cm\noindent [42] Saxton, M.: Anomalous sub-diffusion in fluorescence photobleaching recovery a Monte Carlo Study, {\it Biophys. J.}, {\bf 81}(2001), 2226-2240.

\vskip.3cm\noindent [43] Sokolov, I.M., Chechkin, A.V. and Klafter, J. : Distributed-order fractional kinetics, {\it Acta Phys. Pol.}, {\bf B35}(2004), 1323-1341.

\vskip.3cm\noindent [44] Sokolov, I.M., and Klafter, J. : From diffusion to anomalous diffusion: a century after Einstein's Brownian motion, {\it Chaos}, {\bf 15}(2004), 026103.

\vskip.3cm\noindent [45] Srivastava, H.M. and Daoust, M.C. : A note on the convergence of Kamp\'e de F\'eriet double hypergeometric series, {\it Math. Nachr.}, {\bf 53}(1972), 151-159.

\vskip.3cm\noindent [46] Wilhelmsson, H. and Lazzaro, E. : {\it Reaction-diffusion Problems in the Physics of Hot Plasmas},  Institute of Physics Publishing, Bristol and Philadelphia, 2001.

\vskip.3cm\noindent [47] Wiman A. : Uber den Fundamental satz in der theorie der Functionen $E_{\alpha}(x)$, {\it Acta Math.}, {\bf 29}(1905), 191-201.

\vskip.5cm\noindent{\bf Appendix A: Riemann-Liouville and Riesz-Feller fractional derivatives}

\vskip.3cm The Riemann-Liouville fractional derivative of order $\alpha>0$ is defined as (Samko et al.[32, p.37]; also see, Kilbas et al.[19])
 $${_0D_t^{\alpha}}f(x,t)={{1}\over{\Gamma(n-\alpha)}}{{{\rm d}^n}\over{{\rm d}t^n}}\int_0^t(t-\tau)^{n-\alpha-1}f(x,\tau){\rm d}\tau, ~n=[\alpha]+1, n\in N, t>0\eqno(A1)
 $$where $[\alpha]$ means the integral part of the number $\alpha$.
 \vskip.2cm The Laplace transform of the Riemann-Liouville fractional derivative is given by the following: (see Oldham and Spanier [29, (3.1.3)], Kilbas et al.[20], Podlubny [31], Samko et al.[33], Mainardi [22], Diethlm [4])
 $$\{{_0D_t^{\alpha}}N(x,t);s\}=\cases{s^{\alpha}\tilde{N}(x,s)-\sum_{r=1}^ns^{r-1}{_0D_t^{\alpha-r}}N(x,t)|_{t=0}, n-1<\alpha<n\cr
 {{\partial^n}\over{\partial t^n}}N(x,t), \alpha=n.\cr}\eqno(A2)
 $$This derivative is useful in deriving the solutions of integral equations of fractional order governing certain physical problems of anomalous reaction and anomalous diffusion. In this connection, one can refer to the monograph by Dzherbashyan [5], Podlubny [31], Samko et al.[33], Oldham and Spanier [29], Kilbas et al.[20], Mainardi [22], Diethelm [4] and a recent paper on the subject [28]. Following Feller [8,9], it is convenient to define the Riesz-Feller space-fractional derivative of order $\alpha$ and skewness $\theta$ in terms of its Fourier transform as
 $$\eqalignno{F\{{_xD_{\theta}^{\alpha}}f(x);k\}&=-\psi_{\alpha}^{\theta}(k)f^{*}(k)&(A3)\cr
 \noalign{\hbox{where}}
 \psi_{\alpha}^{\theta}(k)&=|k|^{\alpha}\exp[i(sign~k){{\theta\pi}\over{2}}],~0<\alpha\le 2, |\theta|\le \min\{\alpha,2-\alpha\}.&(A4)\cr}
 $$Further, when $\theta=0$, we have a symmetric operator with respect to $x$ that can be interpreted as
 $${_xD_0^{\alpha}}=-\left[-{{{\rm d}^2}\over{{\rm d}x^2}}\right]^{\alpha/2}.\eqno(A5)
 $$This can be formally deduced by writing $-(k)^{\alpha}=-(k^2)^{\alpha/2}$. For $0<\alpha<2$ and $|\theta|\le\min\{\alpha,2-\alpha\}$, the Riesz-Feller derivative can be shown to possess the following integral representation in $x$ domain:

 $$\eqalignno{{_xD_{\theta}^{\alpha}}f(x)&={{\Gamma(1+\alpha)}\over{\pi}}\bigg\{\sin[(\alpha+\theta)\pi/2]
 \int_0^{\infty}{{f(x+\xi)-f(x)}\over{\xi^{1+\alpha}}}{\rm d}xi\cr
 &+\sin[(\alpha-\theta)\pi/2]\int_0^{\infty}{{f(x-\xi)-f(x)}\over{\xi^{1+\alpha}}}{\rm d}\xi\bigg\}.&(A6)\cr}
 $$For $\theta=0$ the Riesz-Feller fractional derivative becomes the Riesz fractional derivative of order $\alpha$ for $1<\alpha\le 2$ defined by analytic continuation in the whole range $0<\alpha\le 2,\alpha\ne 1$ (see Gorenflo and Mainardi [12]) as
 $$\eqalignno{{_xD_0^{\alpha}}&=-\lambda[I_{+}^{-\alpha}-I_{-}^{-\alpha}]&(A7)\cr
 \noalign{\hbox{where}}
 \lambda&={{1}\over{2\cos(\alpha\pi/2)}}; ~I_{\pm}^{-\alpha}={{{\rm d}^2}\over{{\rm d}x^2}}I_{\pm}^{2-\alpha}.&(A8)\cr}
 $$The Weyl fractional integral operators are defined in the monograph by Samko et al. [33] as

 $$\eqalignno{(I_{+}^{\beta}N)(x)&={{1}\over{\Gamma(\beta)}}\int_{-\infty}^{\infty}(x-\xi)^{\beta-1}N(\xi){\rm d}\xi,~\beta>0\cr
 (I_{-}^{\beta}N)(x)&={{1}\over{\Gamma(\beta)}}\int_x^{\infty}(\xi-x)^{\beta-1}N(\xi){\rm d}\xi,~\beta>0.&(A9)\cr}
 $$
 \vskip.3cm\noindent{\bf Note A1}.\hskip.3cm We note that ${_xD_0^{\alpha}}$ is a pseudo differential operator. In particular, we have
 $$\eqalignno{{_xD_0^2}&={{{\rm d}^2}\over{{\rm d}x^2}}\hbox{  but  }{_xD_0^1}\ne {{{\rm d}}\over{{\rm d}x}}&(A10)\cr
 \noalign{\hbox{For $\theta=0$ we have}}
 F\{{_xD_0^{\alpha}}f(x);k\}&=-|k|^{\alpha}f^{*}(k).&(A11)\cr}
 $$The H-function is defined by means of a Mellin-Barnes type integral in the following manner [25]:
 $$\eqalignno{H_{p,q}^{m,n}(z)&=H_{p,q}^{m,n}\left[x\bigg\vert_{(b_q,B_q)}^{(a_p,A_p)}\right]\cr
 &=H_{p,q}^{m,n}\left[z\bigg\vert_{(b_1,B_1),...,(b_q,B_q)}^{(a_1,A_1),...,(a_p,A_p)}\right]={{1}\over{2\pi i}}\int_{\Omega}\Theta(\xi)z^{-\xi}{\rm d}\xi&(A12)\cr
 \noalign{\hbox{where $i=\sqrt{-1}$ and}}
 \Theta(\xi)&={{\left\{\prod_{j=1}^m\Gamma(b_j+B_j\xi)\right\}\left\{\prod_{j=1}^n\Gamma(1-a_j-A_j\xi)\right\}}
 \over{\left\{\prod_{j=m+1}^q\Gamma(1-b_j-B_j\xi)\right\}\left\{\prod_{j=n+1}^p\Gamma(a_j+A_j\xi)\right\}}}&(A13)\cr}
 $$and an empty product is interpreted as unity; $m,n,p,q\in N_0$ with $0\le n\le p, 1\le m\le q$, $A_i,B_j\in R_{+}$, $a_i,b_j\in C$, $i=1,...,p, j=1,...,q$ such that the poles of $\Gamma(b_j+B_j\xi),j=1,...,m$ are separated from those of $\Gamma(1-a_j-A_j\xi), j=1,...,n$, where $N_0=\{0,1,2,...\}; R=(-\infty,\infty),R_{+}=(0,\infty)$ and $C$ being the complex number field. A comprehensive account of the H-function is available from Mathai et al. [25] and Kilbas et al. [20]. We also need the following result in the analysis that follows: Haubold et al. [14] has shown that
 $$F^{-1}[E_{\beta,\gamma}(-at^{\beta}\psi_{\alpha}^{\theta}(k);x]={{1}\over{\alpha|x|}}
 H_{3,3}^{2,1}\left[{{|x|}\over{(at^{\beta})^{1/\alpha}}}\bigg\vert_{(1,1/\alpha), (1,1), (1,\rho)}^{(1,1/\over\alpha),(\gamma,\beta/\alpha),(1\rho)}\right]\eqno(A14)
 $$where $\Re(\alpha)>0,\Re(\beta)>0,\Re(\gamma)>0$.

 \vskip.3cm\noindent{\bf Appendix B: Convergence of the double power series}

 \vskip.3cm The following lemma given by the authors is needed in the analysis.

 \vskip.3cm\noindent{\bf Lemma}. (Saxena et al. [39]). {\it For all $a,\alpha,\beta>0$ there holds the formula
 $$\sum_{m=0}^{\infty}\sum_{n=0}^{\infty}{{(1)_{m+n}}\over{m!n!}}{{x^my^n}\over{(a)_{\alpha m+\beta n}}}=\Gamma(a)S_{1:0;0}^{1:0;0}\left[x,y\bigg\vert_{[a:\alpha;\beta]:-;-}^{[1:1;1]:-;-}\right]\eqno(B1)
 $$where $S$ stands for the Srivastava-Daoust function (Srivastava et al. [45]), defined by (B2).}

 \vskip.3cm\noindent{\bf Definition B1}.\hskip.3cm Srivastava-Daoust generalization of the Kamp\'e de F\'eriet hypergeometric series in two variables is defined by the double hypergeometric series as [45]:

 $$\eqalignno{&S_{C:D;D'}^{A:B;B'}(x,y)=S_{C:D;D'}^{A:B;B'}\left[
x,y\bigg\vert_{[(c):\delta,\epsilon]:[(d):\eta]:[(d'):\eta']}^{[(a):\theta,\phi]:[(b):\psi]:[(b'):\psi']}\right]\cr
&=\sum_{m=0}^{\infty}\sum_{n=0}^{\infty}{{\{\prod_{j=1}^{A}\Gamma(a_j+m\theta_j+n\phi_j)\}
\{\prod_{j=1}^{B}\Gamma(b_j+m\psi_j)\}\{\prod_{j=1}^{B'}\Gamma(b_j'+n\psi_j')\}x^my^n}\over{\{\prod_{j=1}^{C}\Gamma(c_j+m\delta_j+n\epsilon_j)\}
\{\prod_{j=1}^{D}\Gamma(d_j+m\eta_j)\}\{\prod_{j=1}^{D'}\Gamma(d_j'+n\eta_j')\}m!n!}}&(B2)\cr}
$$where the coefficients $\theta_1,...,\theta_{A},...,\eta'_1,...,\eta'_{D'}>0$. For the sake of brevity $(a)$ is taken to denote the sequence of $A$ parameters $a_1,...,a_{A}$ with a similar interpretation for $(b),...,(d')$. Srivastava and Daoust have shown [45] that the series (B2) converges for all $x,y\in C$ when
$$\eqalignno{\Delta&=1+\sum_{j=1}^{C}\delta_j+\sum_{j=1}^{D}\eta_j-\sum_{j=1}^{A}\theta_j-\sum_{j=1}^{B}\psi_j>0&(B3)\cr
\noalign{\hbox{and}}
\Delta'&=1+\sum_{j=1}^{C}\epsilon_j+\sum_{j=1}^{D'}\eta_j'-\sum_{j=1}^{A}\phi_j-\sum_{j=1}^{B'}\psi_j'>0.&(B4)\cr}
$$For a detailed account of the convergence conditions of the double series, see Srivastava and Daoust [45].

\bye